\documentclass[12pt]{article}

\usepackage[utf8]{inputenc}

\usepackage{amsfonts}
\usepackage{amssymb}
\usepackage{amsmath}
\usepackage{amsthm}
\usepackage[T2A]{fontenc}

\usepackage{amsfonts}
\usepackage{amssymb}
\usepackage{amsmath}
\usepackage{amsthm}

\begin{document} 
\begin{Huge}
\centerline{
Теорема Хинчина о системах Чебышева\footnote{Статья опубликована в журнале
"Историко-математические исследования" вып. 15(50), 55 - 64, 2014.
Здесь мы исправили некоторые опечатки, имеющиеся в публикации.}
}
\vskip+0.5cm
\end{Huge}
\begin{Large}
\centerline{Н.Г. Мощевитин\footnote{Работа поддержана грантами РФФИ №  12-01-00681а и Правительства России, проект 11.
G34.31.0053.}}
\end{Large}
\vskip+2.0cm

В 2006 году были переизданы основные труды А.Я. Хинчина по теории чисел \cite{HHH}.
Большинство из них посвящено теории диофантовых приближений, в частности, многомер- ным диофантовым приближениям.
 Несомненно, что основополагающей в этой области является работа \cite{Paler},
 опубликованная в 1926 году. В ней в зародыше имеется практически вся теория линейных однородных и неоднородных диофантовых приближений, развивав- шаяся впоследствии как самим Хинчиным, так  и другими математиками   -   например, В. Ярником, К. Малером, Дж. Касселсом.
Однако многие классические результаты Хинчина (равно как и его последователей) оказывались забытыми и даже неоднократно передоказывались другими учеными.

В настоящей статье мы сравниваем оригинальный результат Хинчина об эквивалентнос- ти  свойства матрицы быть регулярной и 
свойства матрицы быть матрицей Чебышева (основной результат статьи \cite{HCH}) и его изложение в книге Касселса
"Диофантовы приближе- ния" \cite{cas}.
Оказывается, что формулировка Хинчина и формулировка Касселса различны, хотя и эквивалентны. Это будет  обсуждаться в пунктах 1 и 4 ниже.

Кроме того, мы вносим ясность в  вопрос, как связаны результаты Хинчина 
о теоремах переноса для однородной и неоднородной задач \cite{HCH}, и соответствующие результаты Ярника \cite{J1, J2}.
Об этом  уже написано в нашем обзоре
(пункт 6.2 из \cite{mUMN}). Там сделан акцент на то, что  основной результат Хинчина из работы \cite{HCH} 1948 года
фактически уже имелся в малоизвестной работе Ярника \cite{J1}\footnote{датируемой 1941 годом и изданной в Праге на чешском языке}
и в ее более позднем варианте \cite{J2}. Это действительно так, но, к сожалению, изложение вопроса в \cite{mUMN}
оказалось запутанным и не везде верным\footnote{Например, следствие теоремы 28 (принадлежащей Ярнику)  из пункта  6.2 из \cite{mUMN}
 на самом деле следует не из этой теоремы (как утверждается в \cite{mUMN}), а из другой теоремы Ярника (V\v{e}ta 6 из \cite{J1}),
которая в \cite{mUMN}  оказалась не процитированной. Также теоремы 31 и 32 из \cite{mUMN} не эквивалентны друг другу, как это там  утверждается.}.
В  пункте 2 настоящей работы 
мы дословно приведем формулировки необходимых теорем Ярника и сравним их с результатами Хинчина. Оказывается, что формулировка Ярника ближе к формулировке из Касселса, чем к оригинальной формулировке Хинчина.

{\bf 1.  Теорема о регулярных системах и системах Чебышева: Хинчин и Касселс.}

Пусть  $m$  и  $n$  суть натуральные числа. Положим $ d = m+n$.
Мы	 будем иметь дело с системой из $n$  линейных форм от $m$
целочисленных перенменных.
Для вещественной матрицы
$$\Theta
=\left(
\begin{array}{ccc}
\theta_1^1&\cdots&\theta_1^m\cr
\cdots &\cdots &\cdots \cr
\theta_n^1&\cdots&\theta_n^m
\end{array}
\right)
$$
рассмотрим соотвнтствующую  систему линейных форм
$$
L_j ({\bf x}) =
L_j^\Theta ({\bf x})
=\sum_{i=1}^m\theta^i_j x_i,\,\,\,\, 1\le j \le n.
$$
Здесь  ${\bf x } = (x_1,...,x_m)$  - набор целочисленных переменных.
Через
$^T\!\Theta$  будем обозначать транспонированную матрицу,
соответствующую транспонированную систему линейных форм будем обозначать
$$
^T\!L_i ({\bf y})
=\sum_{j=1}^n\theta^i_j y_j,\,\,\,\, 1\le i \le m.
$$
Целочисленнным набором переменнх будет здесь
${\bf y} = (y_1,...,y_n)$.

Мы будем рассматривать
функции
$$
\psi (t) = \psi^\Theta (t) =
\min_{{\bf x}\in \mathbb{Z}^m:\, 0< \max |x_i|\le t}\,\,\,
\max_{1\le j \le n}
||
L_j({\bf x})||,
$$
$$
^T\!\psi (t) = \psi^{^T\!\Theta } (t) =
\min_{{\bf y}\in \mathbb{Z}^n:\, 0< \max |y_i|\le t}\,\,\,
\max_{1\le i \le m}
||
^T\!L_i({\bf y})||.
$$
Из теоремы Минковского о выпуклом теле следует, что  при $t\ge 1$  выполнено
$$
\psi (t) \le t^{-\frac{m}{n}},\,\,\,\,
^T\!\psi (t) \le t^{-\frac{n}{m}}.
$$

По Хинчину  "{\it если при любом $\varepsilon>0$ и любом достаточно большом $t$ могут быть решены в целых числах неравенства
\begin{equation}\label{pooo}
|L_j| < \frac{1}{t}\,\,\,\,1\le j \le n,\,\,\,
0<\max_{1\le i \le m} |x_i| < \varepsilon t^{\frac{n}{m}},
\end{equation}
то мы называем систему чисел $\theta_j^i$ (или систему линейных форм,
или систему уравнений $L_j = 0$\footnote{добавим, или матрицу $\Theta$}) сингулярной".}

Таким образом, матрица
$\Theta$ является сингулярной, если для любого $\varepsilon '>0$  при достаточно больших $t$  выполнено
$$
\psi (t) < \varepsilon ' t^{-\frac{m}{n}},
$$
или
$$
\limsup_{t\to \infty}  t^{\frac{m}{n}} \psi (t) =0.
$$

Матрицы, не являющиеся сингулярными, называются {\it регулярными}.
Итак, $\Theta$  регулярна, если
при некотором положительном $\varepsilon$  найдется возрастающая к бесконечности последова- тельность $t_k$ такая, что при
$t = t_k$  система неравенств (\ref{pooo}) не имеет решений в целых числах, или, что то же самое,
$$
\limsup_{t\to \infty}  t^{\frac{m}{n}} \psi (t) >0.
$$

Приведем достовно определение Хинчина системы Чебышева \cite{HCH}.
Хинчин пишет\footnote{Мы стараемся, по возможности, приводить точные цитаты; тем не менее, нам
приходится немного 
изменить обозначения, чтобы они соответстввали дург другу всюду в настоящей статье.
Например, в оригинальной работе Хинчина \cite{HCH}
через  $L_j$ обозначалаcь линейная форма
от  $m+1$ переменной
$L_j =\sum_{i=1}^m \theta^i_jx_i-y_j$, и
 неравенства (\ref{0}) 
выглядели следующим образом:
$$
|
L_j-\alpha_j|< \frac{\Gamma}{X^{\frac{m}{n}}}
.$$
}:\\ "{\it Будем называть систему чисел $\theta_j^i$ системой Чебышева, если, каковы бы ни были вещест- веннные числа $\alpha_j \, (1\le j \le n)$,
существует такая положительная постоянная $\Gamma$, что система неравенств
\begin{equation}\label{0}
|
L_j^\Theta ({\bf x}) - y_j-\alpha_j|< \frac{\Gamma}{X^{\frac{m}{n}}}
\end{equation}
имеет целочисленные решения $x_i, y_j, 1\le i \le m, 1\le j \le n$  с каким угодно большим $X=\max_{1\le i \le m } |x_i|$.}"

Итак, по Хинчину система $\Theta$ называется системой Чебышева, если
\begin{equation}\label{1}
\forall \,\alpha \,\,\, \exists \,\Gamma\,\,\, \forall X_0 \,\,\, \exists \,\text{решение системы} \, (\ref{0}) \,\,\,
\text{с} \,\, X \ge X_0. 
\end{equation}

Основным результатом работы Хинчина  \cite{HCH}  является следующая теорема (мы дословно приводим формулировку Хинчина).
 
{\bf Теорема А.} (Хинчин \cite{HCH})\,\, {\it
Регулярность системы чисел $\theta^i_j$  (или системы форм $L_j$, или  системы уравнений $L_j =0$\footnote{или матрицы $\Theta$}) является необходимым и достаточным  условием для того, чтобы эта система была системой Чебышева.}

Дословно приведем формулировку соответствующей теоремы из книги Касселса \cite{cas}.
(теорема XIII из параграфа 7 главы V).
Для этого прежде отметим, что определение сингулярной и регулярной матрицы у Касселса точно такое же, как и приведенное выше определение Хинчина.

{\bf Теорема Б.} (Касселс \cite{cas})\,\, {\it
Для того, чтобы система $L_j({\bf x})$ была регулярна, необходимо и достаточно, чтобы существовало число $\delta>0$,
такое, что неравенство
$$
\left( \max_j ||L_j({\bf x}) -\alpha_j ||\right)^n
\left(
\max_i |x_i|\right)^m <\delta
$$
имело бесконечно много целых решений ${\bf x}$ для каждого действительного $\alpha$.}

Роль величины $\Gamma$ из определения Хинчина у Касслеса играет величина $\delta^{\frac{1}{n}}$.

Мы видим, что по Касселсу, эквивалентным условием регулярности системы  $L_j({\bf x})$
(или матрицы $\Theta$, что то же самое)
является следующее выказываение:
\begin{equation}\label{2}
\exists \,\Gamma \,\,\,
\forall \,\alpha \,\,\,  \forall X_0 \,\,\, \exists \,\text{решение системы} \, (\ref{0}) \,\,\,
\text{с} \,\, X \ge X_0. 
\end{equation}

Неожиданным наблюдением для нас является то, что высказывания (\ref{1}) и (\ref{2}),
вообще-говоря, различные - они отличаются порядком кванторов. 
Как учат студентов, такие высказывания вовсе не обязаны быть равносильными.
Высказывание (\ref{1}) предполагает, что $\Gamma$ зависит от $\alpha$,
в то время как высказывание (\ref{2}) предполагает, что $\Gamma$ одно и то же для всех $\alpha$.

Тем не менее и оригинальное докательство  Хинчина, и доказательство, из книги Касселса являются верными.
Таким образом оказывается, что в рассматриваемой здесь задаче
 (\ref{1}) и (\ref{2}) эквивалентны, и кванторы переставить можно.

Объяснение этого в следующем.
В явном виде можно построить функцию $g:\mathbb{R}_+\to \mathbb{R}_+$
и в более-менее явном виде можно 
указать по матрице $\Theta$   вектор
$\eta^\Theta = (\eta^\Theta_1,...,\eta^\Theta_n) \in \mathbb{R}^n$  такой что  
если  рассмотреть величину
$$
\Delta = \liminf_{{\bf x} \in \mathbb{Z}^m,\,\,|{\bf x}| \to +\infty }
\left( \max_{1\le j\le n} ||L_j({\bf x}) -\eta_j^\Theta ||\right)^n
\left(
\max_i |x_i|\right)^m 
$$
и окажется, что $\Delta <+\infty$,
то для любого вектора
$\alpha = (\alpha_1,...,\alpha_n)\in \mathbb{R}^n$
будет выполнено
$$
\liminf_{{\bf x} \in \mathbb{Z}^m,\,\,|{\bf x}| \to +\infty }
\left( \max_{1\le j\le n} ||L_j({\bf x}) -\alpha_j ||\right)^n
\left(
\max_i |x_i|\right)^m \le g(\Delta ).
$$
В пункте 4 мы приведем точную формулировку этого результата. Но для этого нам потребуется сказать несколько слов о том, как строится множество допустимых векторов $\eta$.
Это мы сделаем в пункте 3.

{\bf 2. Регулярность и свойство Чебышева: Хинчин и Ярник.} 

Для того, чтобы привести в оригинальном виде результат Ярника из работ \cite{J1,J2} нам понадобится функция
$$
\psi^{\Theta, \alpha}(t ) =
\min_{{\bf x}\in \mathbb{Z}^m\setminus\{{\bf 0}\}: \max |x_i| \le t}
\,\,\,\, \max_{1\le j 
\le n} ||L_j({\bf x}) -\alpha_j ||
,$$
определенные для вещественного вектора $\alpha = (\alpha_1,...,\alpha_n)$.

Во всех теоремах ниже предполагается что $\varphi (t) $ и $\rho(t)$ положительные убывающие функции вещественного аргумента, 
причем взаимно обрантые.
Дополнительно предполагается, что
$\varphi (0) = \rho (0) = 0$,
что для некоторого $\eta>0$
функция $ \varphi (t)\cdot t^{-\eta} $ возрастает в интервале $ t>0$ и
$\lim_{t\to \infty}  \varphi (t)\cdot t^{-\eta} =+\infty$.
Кроме того, предполагается что   выполнено
$\varphi (\alpha t) >\alpha^\eta \varphi (t)$ для всех $\alpha > 1, t >0$.

Яринк \cite{J1} доказал следующие две теоремы (V\v{e}ta 5, 7  из \cite{J1}):

{\bf Теорема В.}
\,\,{\it
Если
$$
\limsup_{t\to \infty} \varphi (t)\cdot\,^T\!\psi (t) > a>0,
$$
то
$$
\liminf_{t\to \infty}
\rho (t) 
\,
\sup_{\alpha: \, 0\le \alpha_j \le 1}
\psi ^{\Theta ,\alpha}(t)\le A,
$$
где
$$
A = \frac{3}{2}d! d \max \left( 1,
\left(\frac{d!d}{2a}\right)^{1/\eta}
\right).
$$
}

{\bf Теорема Г.}
\,\,{\it
Если
$$
\limsup_{t\to \infty} \varphi (t)\,\, ^T\!\psi (t)< \infty,
$$
то найдется вектор $\alpha$, такой что
$$
\liminf_{t\to \infty}
\rho (t) \psi^{\Theta,\alpha} (t)  >0.
$$
}
Если положить $\varphi (t) =\gamma t^{\frac{n}{m}}$ то $ \rho (t) =( t/\gamma)^{\frac{m}{n}}, \,\, \gamma >0$,
то происходит следующее.

1. Из теоремы В получаем, что
 {\it 
если матрица  $^T\!\Theta$  регулярна, то матрица  $\Theta$  удовлетворяет высказыванию  (\ref{2}).
}

2. Из теоремы Г получаем, что
{\it если матрица  $^T\!\Theta$  сингулярна, то матрица  $\Theta$  не удовлетворяет высказыванию  (\ref{2}).
}

Итак Ярник в \cite{J1} доказал что регулярность матрицы $^T\!\Theta$ 
эквивалентна высказыванию (\ref{2}) для матрицы $\Theta$. В пункте 6.2 нашей обзорной работы \cite{mUMN} утверждается нечто другое - 
то что Ярник доказал что регулярность матрицы $^T\!\Theta$ эквивалентна
 тому, что $\Theta$  есть матрица Чебышева, то есть удовлетворяет (\ref{1}).

 Последнее, конечно, в каком-то смысле верно, ибо (\ref{1})
и (\ref{2}) действительно эквивалентны.
Но ни у Хинчина, ни у Ярника, ни у Касселса мы не нашли в явном виде утверждения о том, что
(\ref{1}) и (\ref{2})  эквивалентны.

Утверждение о том, что матрица  $\Theta$  сингулярна тогда и только тогда, когда сингулярна матрица  $^T\!\Theta$
в явном виде в работах Ярника \cite{J1,J2} {\it не содержится}.
Это утверждение составляет основной результат работы Хинчина \cite{Spere}  1948 года. Однако оно сразу следует из общей теоремы Малера\footnote{как это и отмечается в главе V, \S 2 книги \cite{cas}}
(Satz 1  из работы \cite{Mcaso} 1939 года).

{\bf 3. Наилучшие приближения и множество   ${\cal N}_\Theta$.}

Для транспонированной матрицы $^T\!\Theta$ 
рассмотрим (конечную или бесконечную)  последо- вательность векторов наилучших приближений
\begin{equation}\label{posl}
{\bf y}_\nu = (y_{\nu,1},...,y_{\nu,n}), \, \nu=1,2,3,... 
\end{equation}
(по поводу определения наилучших приближений и их простейших свойств мы сошлемся на нащу обзорную работу \cite{mUMN}).
Естественно, наиболее интересен и важен случай, когда последовательность (\ref{posl}) бесконечна.

По аналогии с \cite{mUMN} мы будем использовать обозначения
$$
Y_\nu = \max_{1\le j \le n} |y_{\nu, j}|,\,\,\,\,
\zeta_\nu = \max_{1\le i \le m}||\,^T\!L_i ({\bf y}_\nu)||
.$$
Из теоремы Минковского о выпуклом теле сразу следует, что
\begin{equation}\label {mmi}
\zeta_\nu^ mY_{\nu+1}^n \le 1.
\end{equation}

Определим множество
$$
{\cal N}_\Theta =
\{\eta = (\eta_1,...,\eta_n)
\in \mathbb{R}^n:\,\,\,
\inf_{\nu \in \mathbb{Z}_+} ||\eta_1 y_{\nu,1}+...+\eta_ny_{\nu.n}||
>0\}
$$
С одной стороны, это множество достаточно маленькое.
Если последовательность
(\ref{posl}) бесконечна,
то мера Лебега множества ${\cal N}_\Theta$ рана нулю.
С другой стороны это множество не очень маленькое -
оно не пусто, и, более того, размерность Хаусдорфа множества 
${\cal N}_\Theta$  равна размерности всего пространства $\mathbb{R}^n$,
 то есть  $n$.
Можно сказать еще больше,  для произвольной матрицы 
$\Theta$  множество ${\cal N}$  является $1/2$-выигрышным в смысле В.М. Шмидта\footnote{ О выигрышных множествах имеется много литературы;  см.
оригинальную работу В.М. Шмидта \cite{s1} и его книгу \cite{s2}; некоторые  недавние статьи о свойстве выиргышности
процитированны в
\cite{mUMN}  и \cite{marX}}, что есть более сильное свойство, чем обладание полной размерностью Хаусдорфа.

{\bf 4. Новая формулировка.}

Положим
\begin{equation}\label{kaa}
K = \frac{1}{2^{d-1}\sqrt{d}}\, {\rm vol}_{d-1},
\{ {\bf z} = (z_1,...,z_d)\in [-1,1]^d:\,\,\,\,
z_1+...+z_d =0\},
\,\,\,\,\,
\frac{1}{\sqrt{d}}\le K\le \sqrt{\frac{2}{d}},
\end{equation}
\begin{equation}\label{roo}
R_{m,n} = \max_{0<r<1} r^n (1-r)^m = r_0^n(1-r_0)^m \ge \frac{1}{2^d},\,\,\,\, 0<r_0<1
\end{equation}
и определим
$$
G_{m,n} (\varepsilon )=
\frac{d^d(m^mn^n)^{d(d-1)}}{(\varepsilon^d K R_0)^{d(d-1)}}
.
$$

{\bf Теорема 1.}\,\,{\it 
Рассмотрим матрицу $\Theta$  и соответствуюшее ей множество  ${\cal N}_\Theta$.
Пусть $\eta =(\eta_1,...,\eta_n)\in {\cal N}_\Theta$,
и возьмем\footnote{Существует положительная постоянная $\varepsilon_{m,n}$,
такая что для любой матрицы $\Theta$ найдется $\eta \in {\cal N}_\Theta$ с $\varepsilon \ge \varepsilon_{m,n}$ (cм Лемму 2 из главы V в \cite{cas}) .
Это связано с тем, что для любой матрицы $\Theta$ последовательность наилучших приближений в среднем растет лакунарным образом
(см., например, \cite{BL,marX}). Лучшую оценку снизу на $\varepsilon_{m,n}$  можно получить с помощью метода Переса-Шлага (см. оригинальную работу \cite{PS}
или работу \cite{Ro}).  }
 такое положительное $\varepsilon$, что
 $$
\inf_{\nu \in \mathbb{Z}_+} ||\eta_1 y_{\nu,1}+...+\eta_ny_{\nu.n}||\ge \varepsilon.
$$
Рассмотрим величину
$$
\Delta = \Delta (\Theta, \eta)  = \liminf_{{\bf x} \in \mathbb{Z}^m,\,\,|{\bf x}| \to +\infty }
\left( \max_{1\le j\le n} ||L_j({\bf x}) -\eta_j||\right)^n
\left(
\max_i |x_i|\right)^m .
$$
Тогда для любого вектора $\alpha =(\alpha_1,...,\alpha_n) \in \mathbb{R}^n$ выполнено неравенство
$$
\liminf_{{\bf x} \in \mathbb{Z}^m,\,\,|{\bf x}| \to +\infty }
\left( \max_{1\le j\le n} ||L_j({\bf x}) -\alpha_j ||\right)^n
\left(
\max_{1\le i \le m}|x_i|\right)^m \le G_{m,n}(\varepsilon) \Delta ^{d(d-1)}.
$$

}

 В  качестве следствия из теоремы 1 сразу 
получается

{\bf Теорема 2.}\,\,{\it
Для матрицы $\Theta$ выполнено высказывание (\ref{1})
тогда и только тогда, когда
для нее выполнено высказывание (\ref{2}).
}

{\bf 5. Схема доказательства теоремы 1.}

Мы будем следовать классическому доказательству Хинчина \cite{HCH, cas}
Берем $\delta >\Delta$. Из тождества
$$
\eta_1y_1+...+\eta_ny_n =
\sum_{j=1}^n y_j \cdot (\eta_j - L_j ({\bf x})) +
\sum_{i=1}^m x_i \cdot \, ^T\!L_j ({\bf y})
$$
при  ${\bf y} = {\bf y}_\nu$ 
получаем, что
\begin{equation}\label{ceen}
\varepsilon < nY_\nu\max_{1\le j \le n} ||L_j ({\bf x}) - \eta_j||
+ m\zeta_\nu X
.\end{equation}
Давайте в дальнейшем (для определенности) рассматриать только случай, когда последова- тельность  (\ref{posl})
 бесконечна.
Выберем тогда  $\nu$ из условия
$$
(R_0Y_\nu)^{\frac{n}{m}} \le X< (R_0Y_{\nu+1})^{\frac{n}{m}},
$$
где  
$$ R_0 = \frac{n\delta^{\frac{1}{n}}}{\varepsilon r_0},$$
а $r_0$ определено в (\ref{roo}).
Если мы теперь предполагаем, что
$$
\max_{1\le j\le n}
||L_j ({\bf x}) - \eta_j|| \le \left(\frac{\delta}{X^m}\right)^{\frac{1}{n}}
,
$$
то из неравенства (\ref{ceen})  и соотношения  (\ref{mmi})  получаем
\begin{equation}\label{doperenos}
\zeta_\nu^m Y_{\nu+1}^n\ge \omega, 
\end{equation}
где
\begin{equation}\label{amega}
\omega = 
\frac{\varepsilon^{d} r_0^n (1-r_0)^m}{m^mn^n \delta}.
\end{equation}
Это означает, что в $(m+n)$-мерном параллелепипеде
$$
\Pi =
\{
{\bf z} = (x_1,...,x_m,y_1,..,y_n) \in \mathbb{R}^d:\,\,\,
\max_{1\le j \le n} |y_j| < Y_{\nu+1},
\,\,\,
\max_{1\le i \le m}
|^T\!L_i ({\bf y}) - x_i|\le \zeta_\nu\}
$$
нет ненудевых целых точек. 
Согласно принципу переноса\footnote{
Здесь мы используем теоерму К. Малера из \cite{Mcaso}
(излагаемую в главе V  книги \cite{cas}),
точнее, несколько более сильный результат О.Н. Германа \cite{gg},
который дает более точное значение константы. Вообще-говоря в некоторых других
 местах нашего доказательства  возникающие постоянные тоже можно выбирать оптимальнее, чем это сделано у нас.}
 в параллелепипеде
$$^T\!\Pi=
\{
{\bf z} = (x_1,...,x_m,y_1,..,y_n) \in \mathbb{R}^d:\,\,\,
\max_{1\le i \le m} |x_i| < A,
\,\,\,
\max_{1\le j \le n}
|L_i ({\bf x}) - y_j|\le B\}
,$$
где
\begin{equation}\label{aa}
A =K {Y_{\nu+1}^n \zeta_\nu^{m-1}},
\,\,\,\,\,
B = K{Y_{\nu+1}^{n-1}\zeta_\nu^m},
\end{equation}
а $K$ определено в (\ref{kaa}), 
 тоже нет неривиальных целях точек.
Для $d$-мерного объема этого параллелепипеда имеет место равенство
$$
{\rm vol}\, ^T\!\Pi = 2^dA^nB^n = (2K)^d (\zeta_\nu^mY_{\nu+1}^n)^{d-1} \ge (2K)^d \omega^{d-1}.
$$
Рассмотрим  последовательные минимумы $\lambda_1,...,\lambda_{d}$
параллелепипеда  $^T\!\Pi$  по отношению к целочиселнной решетке  $\mathbb{Z}^{d}$. Ясно, что  $\lambda_{d-1}\ge ...\ge \lambda_1 \ge 1$.
Используя теоерму Минковского, видим, что
$$
\lambda_{d} \le \frac{2^{d}}{\lambda_1 ...\lambda_{d-1}{\rm vol}\, ^T\!\Pi} \le \frac{1}{K^d\omega^{d-1}}.
$$
Так как   в параллелепипеде $\lambda_{d}\,\cdot \, ^T\!\Pi$  имеется $d$  независимых целых точек, то параллелепипед\footnote{Скорее всего, множитель   в определении параллелепипеда  $\Pi^*$  выбран не оптимально. Если можно взять меньший множитель, то, естественно,  константа в теореме 1
будет лучше.}
$$
\Pi^* =  \frac{d}{K^d\omega^{d-1}} \cdot\, ^T\!\Pi \supset  d\lambda_{d}\cdot \,^T\!\Pi
$$
cодержит фундаментальную область по отношению к решетке $\mathbb{Z}^{d}$.
Значит для любого вектора  $\xi \in \mathbb{R}^{d}$  трансляция  $\Pi^*+\xi$  содержит некоторую целую точку.
В частности, для любого $\alpha = (\alpha_1,...,\alpha_n) \in \mathbb{R}^n$ найдется целая точка ${\bf x} \in \mathbb{Z}^m$,
такая что
$$
\max_{1\le i \le m} |x_i|\le   \frac{d}{K^d\omega^{d-1}} \cdot A,\,\,\,\,
\max_{1\le j \le n} ||L_j ({\bf x}) - \alpha_j||
\le  \frac{d}{K^d\omega^{d-1}} \cdot B.
$$
Таким образом, с учетом (\ref{mmi},\ref{amega},\ref{aa}) получаем
$$
\left(
\max_{1\le j \le n} ||L_j ({\bf x}) - \alpha_j||
\right)^n\left(\max_{1\le i \le m}|x_i|\right)^m \le
\frac{d^d}{K^{d^2}\omega^{d(d-1)}} A^nB^m =
$$
$$=
\frac{d^d}{(K\omega)^{d(d-1)}}  (\zeta_\nu^m Y_{\nu+1}^n)^{d-1} \le
\frac{d^d}{(K\omega)^{d(d-1)}}   = G_{m,n} \delta^{d(d-1)}.
 $$
Поскольку $\delta>\Delta$ можно выбирать  сколь угодно близким к $\Delta$,
теорема 1  доказана.

\end{document}